\documentstyle[12pt]{amsart}
\title[Macdonald symmetric polynomials]
{A new derivation of the inner product formula
for the Macdonald symmetric polynomials}
\author{
{\sc Katsuhisa Mimachi}
}
\date{}
\numberwithin{equation}{section}
\newtheorem{thm}{Theorem}

\newtheorem{lem}{Lemma}

\begin{document}
\maketitle
\begin{abstract}
We give a short proof of the inner
product conjecture for the symmetric Macdonald polynomials
of type $A_{n-1}$. As a special case, the corresponding
constant term conjecture is also proved.
\end{abstract}
%
%\footnote[0]{
%{\em $1991$ Mathematics Subject Classification\/}.
%Primary ****\,;
%Secondary ****
%33C55, 33D45, 39A10, 39A12
%.}
%
\section{Introduction}
Macdonald's inner product formula, conjectured in
\cite{Mac1}, was recently proved for arbitrary root systems
by Cherednik \cite{C1}, using the double affine Hecke algebras.
In addition to Cherednik's proof, a combinatorial proof by 
Macdonald \cite{Mac1} and representation-theoretic proof 
by Etingof and Kirillov Jr.\cite{EK} have been given for the 
$A_{n-1}$ case. The aim of the present note is to give a short 
proof for the $A_{n-1}$ case by means of asymptotic analysis 
with $q$-Selberg type integrals.  One of our motivations is 
to clarify the argument on the integral representation of 
solutions of eigenvalue problems of the Macdonald type \cite{Mi1}. 
In that case, choice of cycles associated with the integral 
corresponds to the choice of different solutions.  Such study 
on the cycles leads to the present argument, another proof 
of the inner product conjecture for the Macdonald symmetric 
polynomials of type $A_{n-1}$.  Our argument includes a new 
proof of the corresponding constant term conjecture as a 
special case (see also~\cite{Mac2}).\par\medskip
Throughout this note, we consider $q$ as a real number
satisfying $0< q< 1$ and $t=q^k\,,$ where $k\in{\Bbb N}.$
\section{Inner product formula}
We begin recalling some funadamental facts.
For a basic reference, we refer the reader  to
\cite{Mac3}.\par\medskip 
A partition $\lambda$ is a sequence
$\lambda=(\lambda_1\,,\lambda_2\,,\ldots,\lambda_n)$ 
of non-negative integers in
decreasing order; 
$\lambda_1\ge\lambda_2\ge\ldots\ge\lambda_n\ge 0$.  
The number of non-zero elements $\lambda_i$ is called 
the length of $\lambda$, denoted by $l(\lambda)\,.$ 
The sum of the $\lambda_i$ is the weight of $\lambda\,,$
denoted by $|\lambda|$.  Given a partition $\lambda\,,$ we
define the conjugate partition 
$\lambda'=(\lambda_1'\,,\lambda_2'\,,\ldots,\lambda_n')$ by
$\lambda_i'=\text{Card}\{\,j\,;\,\lambda_j\ge i\,\}.$\par  
On partitions, the dominance (or natural) ordering is 
defined by 
$$\lambda\ge \mu \Leftrightarrow |\lambda|=|\mu| \mbox{ and }
\lambda_1+\cdots+\lambda_i\ge\mu_1+\cdots+\mu_i 
\mbox{ for all } 
i\ge 1.$$ \par\medskip
We consider the ring ${\Bbb C}[x]={\Bbb C}[x_1,\ldots,x_n]$
of polynomials in $n$ variables $x=(x_1,\ldots,x_n).$ 
The subring of all symmetric polynomials is denoted by 
 ${\Bbb C}[x]^{\frak S_n}\,.$\par\medskip
For $f=\sum_\beta f_\beta x^\beta\in{\Bbb C}[x]\,,$
we define $$\bar{f}=\sum_\beta f_{\beta} x^{-\beta}$$
and let $[f]_1$ denote the constant term of $f\,.$\par\medskip
The inner product is defined by
$$\langle f, g\rangle=\frac{1}{n!}
[f\bar{g}\Delta]_1$$
for $f,g\in{\Bbb C}[x],$ with 
$$\Delta=\Delta(x)=\prod_{1\le i\ne j \le n}
\frac{(x_i/x_j;q)_\infty}{(tx_i/x_j;q)_\infty}
=\prod_{1\le i\ne j \le n}(x_i/x_j;q)_k\,,$$
where $(a;q)_\infty=\prod_{i\ge 0}(1-aq^i)$ 
and $(a;q)_n=(a;q)_\infty/(q^na;q)_\infty\,.$\par\medskip
Then there is a unique family of symmetric polynomials
$P_\lambda(x)=P_\lambda(x;q,t)\in{\Bbb C}[x]^{\frak S_n}$
indexed by the partition 
$\lambda=(\lambda_1,\ldots,\lambda_n)$ such that
\begin{enumerate}
\item$
P_\lambda=m_{\lambda}+\sum_{\mu<\lambda}
c_{\lambda\mu}m_{\mu}\,,$
\item
$\langle P_\lambda,P_\mu\rangle=0
\quad\text{if}\quad
\lambda\ne \mu\,,
$
\end{enumerate}
where each $m_\mu$ expresses the monomial symmetric polynomial
indexed by $\mu\,.$ The polynomials $P_{\lambda}$ are called 
{\it Macdonald symmetric polynomials} (associated with the root 
system of type $A_{n-1})\,.$\par\medskip
Our aim is to prove the following:
\medskip
\begin{thm}%\label{thm}
We have
\begin{equation*}
\langle P_{\lambda},P_{\lambda}\rangle
=\prod_{1\le i<j\le n}\,\prod_{r=1}^{k-1}
\frac{1-q^{\lambda_i-\lambda_j+r}t^{j-i}}
{1-q^{\lambda_i-\lambda_j-r}t^{j-i}}\,.
\end{equation*}
\end{thm}
\par\medskip
When $\lambda=0$ (so that $P_\lambda=1$), the formula gives
the constant term of $\Delta(x)$. This is the constant term
conjecture of type $A_{n-1}$ (see~\cite{Mac0}).   
\section{Proof of Theorem}
\begin{lem} If $m\ge n,$ for a polynomial 
$\psi(x)=\psi(x_1,\ldots,x_n)$,  we have 
{\allowdisplaybreaks
\begin{align*}
&\hspace{2cm}
\left(\frac{1}{2\pi\sqrt{-1}}\right)^n
\int_{T^n}
\prod
\begin{Sb}
1\le i\le m\\
1\le j\le n
\end{Sb}
\frac{1}{(y_i/x_j;q)_k}\Delta(x) 
\psi(x)\frac{dx_1\cdots dx_n}{x_1\cdots x_n}\\
&=\sum
\begin{Sb}
\{i_1,\ldots,i_n\}\\
\subset \{1,\ldots,n\}
\end{Sb}
\quad
\sum_{0\le l_1,\cdots,l_n\le k-1}\\
& \hspace{1cm}
\text{ Res }\begin{Sb}
x=(y_{i_1}q^{l_1},\ldots,y_{i_n}q^{l_n})
\end{Sb}
\left\{
\prod
\begin{Sb}
1\le i\le m\\
1\le j\le n
\end{Sb}
\frac{1}{(y_i/x_j;q)_k}\Delta(x) 
\psi(x)\frac{dx_1\cdots dx_n}{x_1\cdots x_n}
\right\}\,,
\end{align*}
}\noindent
where $i_1,\ldots,i_n$ are distinct, and
$T^n=\{(t_1,\ldots,t_n)\in {\Bbb C}^n ; |t_i|=1  (1\le i\le n)\}\,.$ 
\end{lem}
\begin{pf}
For a polynomial $\psi(x_1,x_2)$ and $0\le l\le k-1,$ we have the equality
\begin{equation}
\begin{split}
&\text{Res}_{x_1=yq^l}
\frac{(x_1/x_2;q)_k(x_2/x_1;q)_k}{(y/x_1;q)_k(y/x_2;q)_k}
\psi (x_1,x_2)\frac{dx_1}{x_1}\frac{dx_2}{x_2}\\
&=
\frac{({yq^l}/{x_2};q)_k({x_2q^{-l}}/y;q)_k}
{(q^{-l};q)_l(q;q)_{k-1-l}(y/x_2;q)_k}
\psi (yq^l,x_2) \frac{dx_2}{x_2}\,.
\end{split}\label{eq3.1}
\end{equation}
Because $(y/x_2;q)_k$ divides 
$({yq^l}/{x_2};q)_k({x_2q^{-l}}/{y};q)_k\,,$
the 1-form (\ref{eq3.1}) has no poles on the $x_2-$plane.
This shows that the set of poles of 
$$\prod
\begin{Sb}
1\le i\le m\\
1\le j\le 2
\end{Sb}
\frac{1}{(y_i/x_j;q)_k}\Delta(x_1,x_2) 
\psi(x_1,x_2)\frac{dx_1dx_2}{x_1x_2}
$$
is the union of $(x_1,x_2)=(y_{i_1}q^{l},y_{i_2}q^{l})$
for $1\le i_1\ne i_2\le m$ and $0\le l\le k-1\,,$ 
which implies the assertion of the above Lemma in the $n=2$ case.  
Repeating this procedure, we have the desired result in case of
general $n.$ 
\end{pf}
It is known ((3.11) in \cite{Mac1}) that
\begin{equation}
\sum_{\lambda}b_{\lambda}
P_{\lambda}(y)P_{\lambda}(x)
=\prod
\begin{Sb}
1\le i\le m\\
1\le j\le n
\end{Sb}
\frac{(ty_{i}x_{j};q)_\infty}{
(y_{i}x_{j};q)_\infty}
=\prod
\begin{Sb}
1\le i\le m\\
1\le j\le n
\end{Sb}
\frac{1}{(y_{i}x_{j};q)_k}\label{eq3.2}
\end{equation}
with
\begin{eqnarray*}
b_{\lambda}=b_{\lambda}(q,t)=\prod_{s\in\lambda}
\frac{1-q^{a(s)}t^{l(s)+1} }{1-q^{a(s)+1}t^{l(s)} }\,.
\end{eqnarray*}
Here the sum is taken over all partitions $\lambda$ 
such that $l(\lambda)\le\text{min}\{m,n\}\,,$
and the arm-length $a(s)$ (resp.
the leg-length $l(s)$) is defined by
$a(s)=\lambda_i -j$ (resp. $l(s)=\lambda_j' -i$ ) 
for a square 
$s=(i,j)$ in the diagram $\lambda\,.$ \par\medskip
The formula (\ref{eq3.2}) in the $m=n$ case with the orthogonality relation gives
{\allowdisplaybreaks
\begin{align}
&b_{\lambda}P_{\lambda}(y)
\langle P_{\lambda},P_{\lambda}\rangle
=\frac{1}{n!}\left(\frac{1}{2\pi\sqrt{-1}}\right)^n
\int_{T^n}\prod_{1\le i,j\le n}
\frac{1}{(y_i/x_j;q)_k}
P_{\lambda}(x)\Delta(x) 
\frac{dx_1\cdots dx_n}{x_1\cdots x_n}\label{eq3.3} \\\nonumber
&\hspace{3mm}=\frac{1}{n!}\sum_{\sigma\in{\frak S}_n}
\quad\sum_{0\le l_1,\cdots,l_n\le k-1}\\\nonumber
&\hspace{1cm}\text{ Res }\begin{Sb}
x=(y_{\sigma(1)}q^{l_1},\ldots,y_{\sigma(n)}q^{l_n})
\end{Sb}
\left\{\prod_{1\le i,j\le n}
\frac{1}{(y_i/x_j;q)_k} 
P_{\lambda}(x)\Delta(x)\frac{dx_1\cdots dx_n}{x_1\cdots x_n}\right\}\\\nonumber
&\hspace{3mm}=\sum_{0\le l_1,\cdots,l_n\le k-1}\\\nonumber
&\hspace{1cm}
\text{ Res }\begin{Sb}
x=(y_1q^{l_1},\ldots,y_nq^{l_n})
\end{Sb}
\left\{
\prod_{1\le i,j\le n}\frac{1}{(y_i/x_j;q)_k}
P_{\lambda}(x)\Delta(x) 
\frac{dx_1\cdots dx_n}{x_1\cdots x_n}\right\}\,.
\end{align}}\noindent
Here the second equality is given by Lemma above 
and the third equality by the symmetry of the
summand with respect to the variables $x=(x_1,\ldots,x_n).$ \par
Next, by changing the integration variables on
the right-hand side according to $x_i\to y_ix_i\,,$ we have
{\allowdisplaybreaks
\begin{align*}
&\sum_{0\le l_1,\cdots,l_n\le k-1}
\text{ Res }
\begin{Sb}
x=(q^{l_1},\ldots,q^{l_n})
\end{Sb}
\Biggl\{\prod_{1\le i,j\le n}
\frac{1}{\displaystyle\biggl(\frac{y_i}{y_jx_j};q\biggr)_k}
\prod_{1\le i\ne j\le n}
{\displaystyle\biggl(\frac{y_ix_i}{y_jx_j};q\biggr)_k} \\
&\hspace{5mm}
P_{\lambda}(y_1x_1,\ldots,y_nx_n)\,
\frac{dx_1\cdots dx_n}{x_1\cdots x_n}  \Biggr\}\,,
\end{align*}}
which tends to
\begin{align*}
&\sum_{0\le l_1,\cdots,l_n\le k-1}
\text{ Res }
\begin{Sb}
x=(q^{l_1},\ldots,q^{l_n})
\end{Sb}\Biggl\{
\frac{x_1^k (x_1x_2)^k\cdots (x_1\cdots x_{n-1})^k}
{\prod_{i=1}^{n}(1/x_i)_k}\\
&\hspace{10mm}
\{\,(y_1x_1)^{\lambda_1}\cdots (y_nx_n)^{\lambda_n}
+\text{lower order terms}\,\}
\frac{dx_1\cdots dx_n}{x_1\cdots x_n}\Biggr\}\\
&=
y^{\lambda}
\sum_{0\le l_1,\cdots,l_n\le k-1}
\text{ Res }
\begin{Sb}
x=(q^{l_1},\ldots,q^{l_n})
\end{Sb}
\Biggl\{
\prod_{i=1}^n
\frac{(x_i)^{\lambda_i+(n-i)k}}{(1/x_i;q)_k}
\frac{dx_1\cdots dx_n}{x_1\cdots x_n}\Biggr\}\\[10pt]
&\hspace{5mm}+\text{lower order terms}\\[10pt]
&=y^{\lambda}
\prod_{i=1}^{n}
\frac{(q^{\lambda_i+(n-i)k+1};q)_{k-1}}{(q;q)_{k-1}}
+\text{lower order terms}\,,
\end{align*}
if
\begin{equation}
1>|y_1|\gg |y_2|\gg\cdots\gg |y_n|\,.\label{eq3.4}
\end{equation}
Comparing the coefficients of $y^\lambda$ of (\ref{eq3.3})
in the region (\ref{eq3.4}) leads to
$$b_\lambda
\langle P_{\lambda},P_{\lambda}\rangle
=\prod_{i=1}^{n}
\frac{(q^{\lambda_i+(n-i)k+1};q)_{k-1}}{(q;q)_{k-1}}\,,
$$
which is equivalent to
$$
\langle P_{\lambda},P_{\lambda}\rangle
=\prod_{1\le i< j\le n}
\frac{(q^{\lambda_i-\lambda_j+1+(j-i)k};q)_{k-1}}
{(q^{\lambda_i-\lambda_j+1+(j-i-1)k};q)_{k-1}}.
$$
Here we used the equality 
$$
b_\lambda=\prod_{1\le i< j\le n}
\frac{(q^{\lambda_i-\lambda_j+1+(j-i-1)k};q)_{k-1}}
{(q^{\lambda_i-\lambda_j+1+(j-i)k};q)_{k-1}}
\prod_{i=1}^n
\frac{(q^{\lambda_i+1+k(n-i)};q)_{k-1}}{(q;q)_{k-1}}\,.
$$
This completes the proof of our Theorem.\par\medskip
{\it Remark.} When we would like to consider the $q=1$ case
directly, we need only modify the proof of Lemma and
the calculation of the residue at the final step.\par\bigskip
{\it Acknowledgement.} The author wishes to thank Professor Masatoshi
Noumi for suggesting the link with the inner product formula after 
the completion of the work in \cite{Mi2}.

\bigskip
\begin{flushleft}
\begin{sc}
Katsuhisa Mimachi\\
Department of Mathematics\\
Kyushu University\\
Hakozaki 33, Fukuoka 812\\
Japan\\
\end{sc}
\smallskip 
{\it E-mail address\/}: mimachi{\char'100}
math.kyushu-u.ac.jp\\
\end{flushleft}
\end{document}